\documentclass{article}
\usepackage{graphics, graphicx, latexsym, amssymb}
\newtheorem{theo}{Theorem}

\begin{document}

\title{Polyominoes with nearly convex columns: \\
An undirected model}
\author{Svjetlan Fereti\'{c} \footnote{e-mail: svjetlan.feretic@gradri.hr} \\ 
Faculty of Civil Engineering, University of Rijeka, \\ 
Viktora Cara Emina 5, 51000 Rijeka, Croatia
\and
Anthony J. Guttmann \footnote{e-mail: tonyg@ms.unimelb.edu.au} \\
ARC Centre of Excellence for\\
Mathematics and Statistics of Complex Systems, \\
Department of Mathematics and Statistics, \\
The University of Melbourne, \\
Parkville, Victoria 3010, Australia}
\maketitle

\begin{abstract}
Column-convex polyominoes were introduced in 1950's by Temperley, a mathematical physicist working on ``lattice gases". By now, column-convex polyominoes are a popular and well-understood model. There exist several generalizations of column-convex polyominoes; an example is a model called \textit{multi-directed animals}.
In this paper, we introduce a new sequence of supersets of column-convex polyominoes. Our model (we call it \textit{level} $m$ \textit{column-subconvex polyominoes}) is defined in a simple way. We focus on the case when cells are hexagons and we compute the area generating functions for the levels one and two. Both of those generating functions are complicated $q$-series, whereas the area generating function of column-convex polyominoes is a rational function. The growth constants of level one and level two column-subconvex polyominoes are $4.319139$ and $4.509480$, respectively. For comparison, the growth constants of column-convex polyominoes, multi-directed animals and all polyominoes are $3.863131$, $4.587894$ and $5.183148$, respectively.
\end{abstract}

\vspace{8mm}

\vspace{3mm}

\vspace{3mm}
\noindent \textit{Suggested running head:} Column-subconvex polyominoes

\section{Introduction}

The enumeration of polyominoes is a topic of great interest to chemists, physicists and combinatorialists alike. In chemical terms, any polyomino (with hexagonal cells) is a possible benzenoid hydrocarbon. In physics, determining the number of $n$-celled polyominoes is related to the study of two-dimensional percolation phenomena. In combinatorics, polyominoes are of interest in their own right because several polyomino models have good-looking exact solutions.

Known results for polyominoes include the fact that the number of $n$-celled polyominoes grows exponentially. More precisely, if $a_n$ denotes the number of $n$-celled polyominoes, then
\begin{itemize}
\item
$\lim_{n \to \infty} a_n^{1/n} = \tau = \sup_{n \ge 1} a_n^{1/n}$
\item
$\lim_{n \to \infty} a_{n+1}/a_n = \tau.$
\end{itemize}
The first result follows from standard concatenation arguments, see e.g. \cite{Voege, book}, while the second result, due to Madras \cite{Madras} relies on a pattern theorem for lattice animals. These results are quite general, and apply {\em mutatis mutandis} to the new polyomino models we consider here.

One can also obtain rigorous bounds on the growth constant $\tau.$ For example, for hexagonal polyominoes we have \cite{Voege}
$$4.8049 \le \tau \approx 5.183148 \le 5.9047.$$ A lower bound is immediately obtainable from the first itemised equation above, and it can be improved with rather more work. The upper bound is obtained by a method due to Klarner and Rivest \cite{KR73}, which relies on mapping each polyomino onto a tree on the dual lattice, and relaxing the rules for tree construction so that over-counting results.

One very popular polyomino model is that of \textit{column-convex polyominoes}. Column-convex polyominoes with hexagonal cells have a rational area generating function. That generating function was found by Klarner in 1967 \cite{Klarner}. The \textit{growth constant} of hexagonal-celled column-convex polyominoes is $3.863131$. (By the growth constant we mean the limit $\lim_{n\rightarrow\infty} \sqrt[n]{a_n}$, where $a_n$ denotes the number of $n$-celled elements in a given set of polyominoes.)

In a previous paper \cite{semi}, one of us (Fereti\'{c}) began to search for polyomino models which are more general than column-convex polyominoes, but still have reasonably simple area generating functions. In \cite{semi}, Fereti\'{c} introduced \textit{level} $m$ \textit{cheesy polyominoes} ($m=1,\: 2,\: 3,\ldots$), and here we shall introduce another sequence of models, which we call \textit{level} $m$ \textit{column-subconvex polyominoes} ($m=1,\: 2,\: 3,\ldots$)\footnote{Bousquet-M\'{e}lou and Rechnitzer's \textit{multi-directed animals} \cite{Rechnitzer} are also a superset of column-convex polyominoes with hexagonal cells.}.

At every level, cheesy polyominoes have a rational area generating function, whereas column-subconvex polyominoes have an area generating function which is unlikely to be algebraic, and indeed, unlikely to be differentiably finite \cite{JPhysA}. Further, at any given level, cheesy polyominoes are an exponentially small subset of column-subconvex polyominoes. The latter set of polyominoes has a greater growth constant than the former set. For example, the growth constant of level one cheesy polyominoes is $4.114908$, while the growth constant of level one column-subconvex polyominoes is $4.319139$. In addition, if we reflect a column-subconvex polyomino about a vertical axis, we get a column-subconvex polyomino again. This kind of invariance under reflection is enjoyed by column-convex polyominoes, but not by cheesy polyominoes. Admittedly, counting level $m$ column-subconvex polyominoes requires more effort than counting level $m$ cheesy polyominoes. Anyway, at level one, column-subconvex polyominoes are not very hard to count. Just as with cheesy polyominoes, as the level increases, the computations quickly increase in size. 

In this paper, the level one column-subconvex model is solved in full detail. We also solved the level two column-subconvex model. To see the level two result (stated with no proof), the reader may visit the web page \cite{worksheet}. The said result involves too many auxiliary expressions to be stated in this paper. (To be specific, there are $33$ auxiliary expressions, of which $25$ are polynomials; the degree of those polynomials is between $20$ and $23$.)

Our computations are done by using Bousquet-M\'{e}lou's \cite{Bousquet} and Svrtan's~\cite{Svrtan} ``turbo'' version of the Temperley method \cite{Temperley}. 

If the reader would like to have more information on the history of polyomino enumeration, or on the role which polyominoes play in physics and chemistry, then he/she may refer to Bousquet-M\'{e}lou's habilitation thesis \cite{habilitation}, or to the book \cite{book}.

\section{Definitions and conventions}

There are three regular tilings of the Euclidean plane, namely the triangular tiling, the square tiling, and the hexagonal tiling. We adopt the convention that every square tile or hexagonal tile has two horizontal edges. In a regular tiling, a tile is often referred to as a \textit{cell}. A plane figure $P$ is a polyomino if $P$ is a union of finitely many cells and the interior of $P$ is connected. See Figure~1. Observe that, if a union of \textit{hexagonal} cells is connected, then it possesses a connected interior as well. 

\begin{figure}
\begin{center}
\includegraphics[width=55mm]{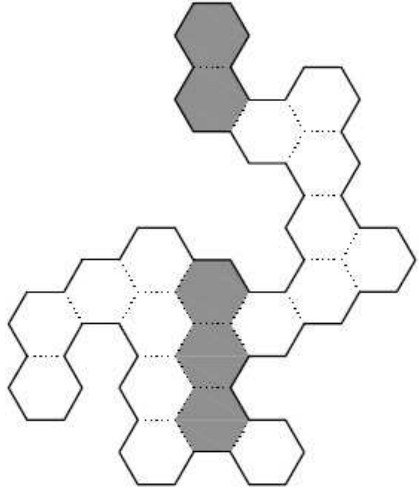}
\caption{A hexagonal-celled polyomino.}
\end{center}
\end{figure}

Let $P$ and $Q$ be two polyominoes. We consider $P$ and $Q$ to be equal if and only if there exists a translation $f$ such that $f(P)=Q$.

If a polyomino $P$ is made up of $n$ cells, we say that the \textit{area} of $P$ is $n$.

Let $R$ be a set of polyominoes. By the \textit{area generating function} of $R$ we mean the formal sum

\begin{displaymath}
\sum_{P \in R} q^{area \ of \ P}.
\end{displaymath}

From now on, we concentrate on the hexagonal tiling. When we write ``a polyomino'', we actually mean ``a hexagonal-celled polyomino''.

Given a polyomino $P$, it is useful to partition the cells of $P$ according to their horizontal projection. Each block of that partition is a \textit{column} of $P$. Note that a column of a polyomino is not necessarily a connected set. An example of this is the highlighted column in Figure 1. On the other hand, it may happen that every column of a polyomino $P$ is a connected set. In this case, the polyomino $P$ is a \textit{column-convex polyomino}. See Figure 2.

\begin{figure}
\begin{center}
\includegraphics[width=47.5mm]{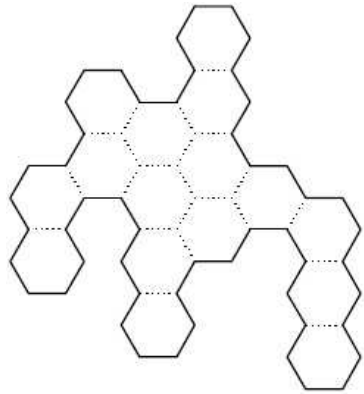}
\caption{A column-convex polyomino.}
\end{center}
\end{figure}

Let $a$ be a column of a polyomino $P$. By the \textit{height} of $a$ we mean the number of those cells which make up $a$ plus the number of those (zero or more) cells which make up the gaps of $a$. For example, in Figure 1, the highlighted column has height $7$, and the next column to the left has height $4$. 

A finite union of cells $P$ is a \textit{level} $m$ \textit{column-subconvex polyomino} if the following holds:

\begin{itemize}
\item $P$ is a polyomino,
\item every column of $P$ has at most two connected components,
\item if a column of $P$ has two connected components, then the gap between the components consists of at most $m$ cells.
\end{itemize}

See Figure 3.

\begin{figure}
\begin{center}
\includegraphics[width=55mm]{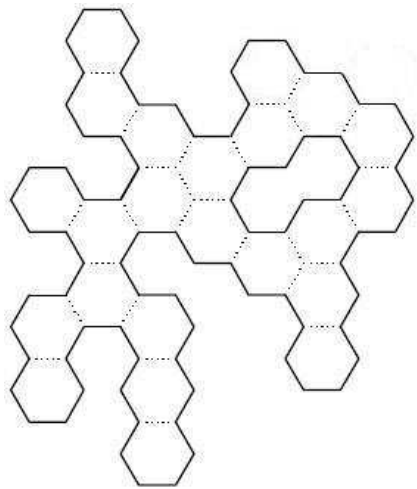}
\caption{A level one column-subconvex polyomino.}
\end{center}
\end{figure}

Let $S$ denote the set of all level one column-subconvex polyominoes.

Let $P$ be an element of $S$ and let $P$ have at least two columns. Then we define the \textit{pivot cell} of $P$ to be the lower right neighbour of the lowest cell of the second last column of $P$. See Figure 4. Observe that the pivot cell of $P$ is not necessarily contained in $P$. 

\begin{figure}
\begin{center}
\includegraphics[width=\textwidth]{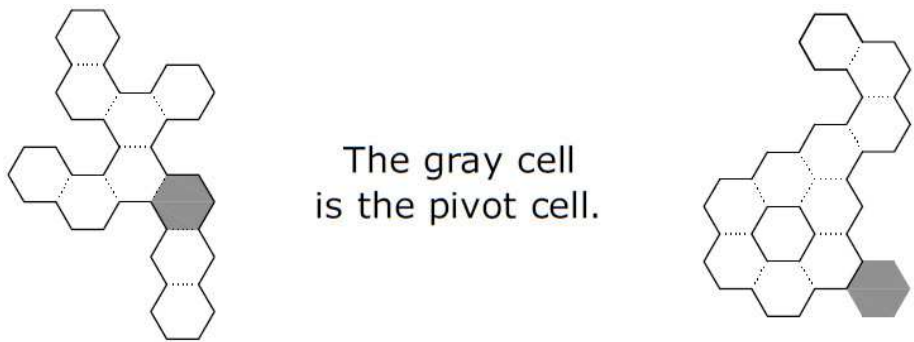}
\caption{The pivot cell.}
\end{center}
\end{figure}

When we build a column-convex polyomino from left to right, adding one column at a time, every intermediate figure is a column-convex polyomino itself. However, when we build a column-subconvex polyomino, this is no longer the case. A ``left factor'' of an element of $S$ need not itself be a polyomino, and therefore need not be an element of $S$.

We say that a figure $P$ is an \textit{incomplete level one column-subconvex polyomino} if $P$ itself is not an element of $S$, but $P$ is a ``left factor'' of an element of $S$. Notice that, if $P$ is an incomplete level one column-subconvex polyomino, then the last (\textit{i.e.}, the rightmost) column of $P$ necessarily has a hole.

Let $T$ denote the set of all incomplete level one column-subconvex polyominoes.

Let $P$ be an element of $S \cup T$ and let $P$ have at least two columns. Then we define the \textit{body} of $P$ to be all of $P$, except the rightmost column of $P$.

Let $P$ be an element of $T$ and let $P$ have at least two columns. We define the \textit{lower pivot cell} of $P$ to be the lower right neighbour of the lowest cell of the second last column of $P$. In addition, we define the \textit{upper pivot cell} of $P$ to be the upper right neighbour of the highest cell of the second last column of $P$.

\section{Notations for generating functions. Partitions of the sets \newline $S$ and $T$}

We shall deal with the following generating functions:

\begin{displaymath}
A(q,t)=\sum_{P \in S} q^{area\ of\ P} \cdot t^{the\ height\ of\ the\ last\ column\ of\ P},
\end{displaymath}

\begin{displaymath}
A_1=A(q,1), \qquad B_1=\frac{\partial A}{\partial t}(q,1),
\end{displaymath}

\begin{eqnarray*}
C(q,u,v) & = & \sum_{P \in T} q^{area\ of\ P} \cdot u^{{the\ height\ of\ the \ upper \atop component\ of\ the\ last\ column\ of\ P} \atop } \cdot \\
& & \cdot v^{{the\ height\ of\ the \ lower \atop component\ of\ the\ last\ column\ of\ P} \atop },
\end{eqnarray*}

\begin{displaymath}
D(u)=C(q,u,1), \qquad E(v)=C(q,1,v), \qquad C_1=C(q,1,1).
\end{displaymath}

Functional equations for the generating functions will be obtained by the ``divide and conquer'' strategy. Namely, now we are going to partition the sets $S$ and $T$.

Let $S_{\alpha}$ be the set of level one column-subconvex polyominoes which have only one column.

Let

\begin{eqnarray*}
S_{\beta} & = & \{P \in S \setminus S_{\alpha}: \mathrm{the \ body \ of \ } P \mathrm{\ lies \ in \ } S\mathrm{, \ the \ last \ column \ of \ } P \\
& & \mathrm{\ has \ no \ hole, \ and \ the \ pivot \ cell \ of \ } P \mathrm{\ is \ contained \ in \ } P \}, \\
S_{\gamma} & = & \{P \in S \setminus S_{\alpha}: \mathrm{the \ body \ of \ } P \mathrm{\ lies \ in \ } S\mathrm{, \ the \ last \ column \ of \ } P \\
& & \mathrm{\ has \ no \ hole, \ and \ the \ pivot \ cell \ of \ } P \mathrm{\ is \ not \ contained \ in \ } P \}, \\
S_{\delta} & = & \{P \in S \setminus S_{\alpha}: \mathrm{the \ body \ of \ } P \mathrm{\ lies \ in \ } S\mathrm{, \ and \ the \ last \ column \ of \ } P \\
& & \mathrm{\ has \ a \ hole} \}, \\
S_{\epsilon} & = & \{P \in S \setminus S_{\alpha}: \mathrm{the \ body \ of \ } P \mathrm{\ lies \ in \ } T\mathrm{, \ and \ the \ last \ column \ of \ } P \\
& & \mathrm{\ has \ no \ hole\} \quad and} \\
S_{\zeta} & = & \{P \in S \setminus S_{\alpha}: \mathrm{the \ body \ of \ } P \mathrm{\ lies \ in \ } T\mathrm{, \ and \ the \ last \ column \ of \ } P \\
& & \mathrm{\ has \ a \ hole} \}.
\end{eqnarray*}

The sets $S_{\alpha}$, $S_{\beta}$, $S_{\gamma}$, $S_{\delta}$, $S_{\epsilon}$ and $S_{\zeta}$ form a partition of $S$. We write $A_{\alpha}$, $A_{\beta}$, $A_{\gamma}$, $A_{\delta}$, $A_{\epsilon}$ and $A_{\zeta}$ for the parts of the series $A$ that come from the sets $S_{\alpha}$, $S_{\beta}$, $S_{\gamma}$, $S_{\delta}$, $S_{\epsilon}$ and $S_{\zeta}$, respectively.

We proceed to the set $T$. We write $T_{\alpha}$ for the set of incomplete level one column-subconvex polyominoes which have only one column. Let $P \in T \setminus T_{\alpha}$. If the body of $P$ lies in $S$, then the said body is in contact with just one of the two connected components of $P$'s last column. The non-contacting component of the last column is located either wholly above or wholly below the second last column of $P$. Let

\begin{eqnarray*}
T_{\beta} & = & \{P \in T \setminus T_{\alpha}: \mathrm{the \ body \ of \ } P \mathrm{\ lies \ in \ } S\mathrm{, \ and \ the \ hole \ of \ the \ last} \\
& & \mathrm{column \ of \ } P \mathrm{\ coincides \ either \ with \ the \ lower \ pivot \ cell \ of \ } P \\
& & \mathrm{or \ with \ the \ upper \ pivot \ cell \ of \ } P \} \quad \mathrm{and} \\
T_{\gamma} & = & \{P \in T \setminus T_{\alpha}: \mathrm{the \ body \ of \ } P \mathrm{\ lies \ in \ } S\mathrm{, \ and \ the \ hole \ of \ the \ last} \\
& & \mathrm{column \ of \ } P \mathrm{\ lies \ either \ below \ the \ lower \ pivot \ cell \ of \ } P \\
& & \mathrm{or \ above \ the \ upper \ pivot \ cell \ of \ } P \}.
\end{eqnarray*}

Let us move on to the case when the body of $P \in T \setminus T_{\alpha}$ lies in $T$. Then the second last column of $P$ has two connected components. It is easy to see that each of those two components must be in contact with the last column of $P$. (This does not mean that each of the two connected components of the last column of $P$ must be in contact with the second last column of $P$.) Now, it may or may not happen that one connected component of $P$'s last column is in contact with both connected components of $P$'s second last column. Accordingly, we define the following two sets:

\begin{eqnarray*}
T_{\delta} & = & \{P \in T \setminus T_{\alpha}: \mathrm{the \ body \ of \ } P \mathrm{\ lies \ in \ } T\mathrm{, \ and \ the \ hole \ of \ the \ last} \\
& & \mathrm{column \ of \ } P \mathrm{\ touches \ the \ hole \ of \ the \ second \ last \ column \ of \ } P \} \quad \mathrm{and} \\
T_{\epsilon} & = & \{P \in T \setminus T_{\alpha}: \mathrm{the \ body \ of \ } P \mathrm{\ lies \ in \ } T\mathrm{, \ and \ the \ hole \ of \ the \ last} \\
& & \mathrm{column \ of \ } P \mathrm{\ does \ not \ touch \ the \ hole \ of \ the \ second \ last \ column \ of \ } P \}.
\end{eqnarray*}

The sets $T_{\alpha}$, $T_{\beta}$, $T_{\gamma}$, $T_{\delta}$ and $T_{\epsilon}$ form a partition of $T$. We write $C_{\alpha}$, $C_{\beta}$, $C_{\gamma}$, $C_{\delta}$ and $C_{\epsilon}$ for the parts of the series $C$ that come from the sets $T_{\alpha}$, $T_{\beta}$, $T_{\gamma}$, $T_{\delta}$ and $T_{\epsilon}$, respectively.

\section{Setting up the functional equations for $A$, $A_1$ and $B$}

To begin with, it is clear that

\begin{equation}
A_{\alpha}=qt+(qt)^2+(qt)^3+\ldots =\frac{qt}{1-qt}.
\end{equation}

If a polyomino $P$ lies in $S_{\beta}$, then the last column of $P$ is made up of the pivot cell, of $i \in \{0,\: 1,\: 2,\: 3,\ldots \: \}$ cells lying below the pivot cell, and of $j \in \{0,\: 1,\: 2,\: 3,\ldots \: \}$ cells lying above the pivot cell. See Figure 5. Hence,

\begin{equation}
A_{\beta} = A_1 \cdot qt \cdot \left[ \sum_{i=0}^{\infty} (qt)^i \right] \cdot \left[ \sum_{j=0}^{\infty} (qt)^j \right] = \frac{qt}{(1-qt)^2} \cdot A_1.
\end{equation}

\begin{figure}
\begin{center}
\includegraphics[width=114mm]{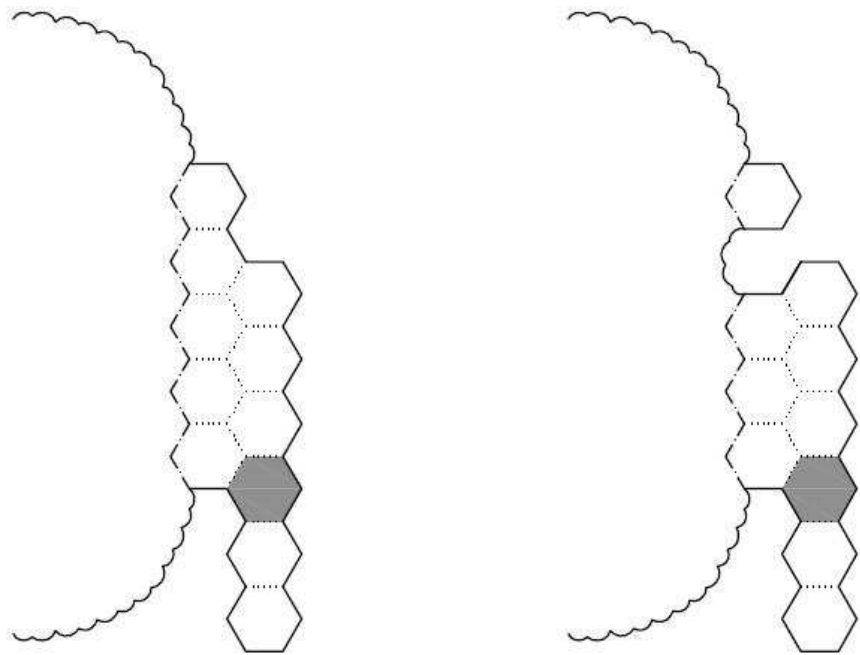}
\caption{The last two columns of two elements of $S_{\beta}$.}
\end{center}
\end{figure}

Consider the following situation. A polyomino $P \in S$ ends with a column $I$. We are creating a new column to the right of $I$, and the result should be an element of $S_{\gamma}$. Then, whether or not the column $I$ has a hole, we can put the lowest cell of the new column in exactly $m$ places, where $m$ is the height of $I$. See Figure 6. Hence

\begin{equation}
A_{\gamma}=\frac{qt}{1-qt} \cdot B_1.
\end{equation}

\begin{figure}
\begin{center}
\includegraphics[width=114mm]{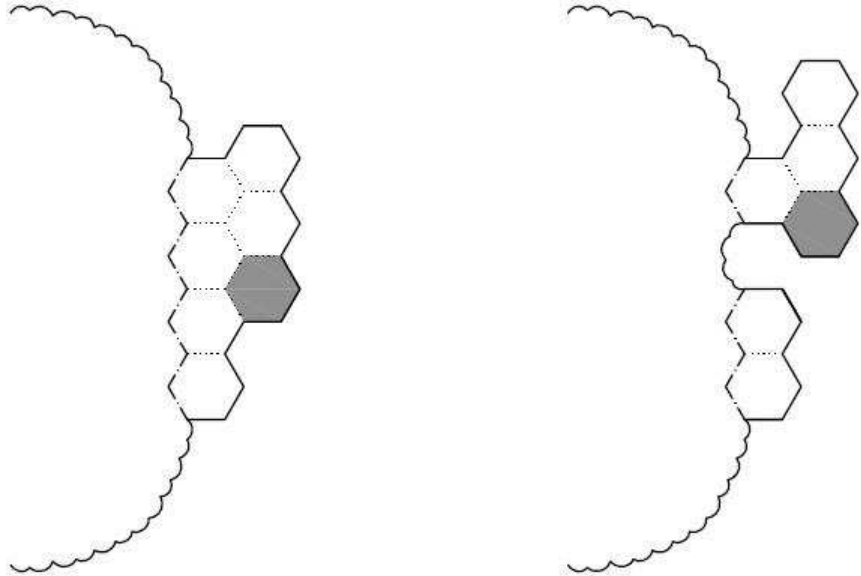}
\caption{The last two columns of two elements of $S_{\gamma}$.}
\end{center}
\end{figure}
 
Let us proceed to another situation. A polyomino $P \in S$ ends with a column $J$. We are creating a new column to the right of $J$, and the result should be an element of $S_{\delta}$. Then, whether or not the column $J$ has a hole, we can put the hole of the new column in exactly $n-1$ places, where $n$ is the height of $J$. See Figure 7. The new column is made up of $i \in \{1,\: 2,\: 3,\ldots \: \}$ cells lying below the hole, of a hole of height one, and of $j \in \{1,\: 2,\: 3,\ldots \: \}$ cells lying above the hole. Altogether,

\begin{equation}
A_{\delta}=\frac{qt}{1-qt} \cdot t \cdot \frac{qt}{1-qt} \cdot (B_1-A_1) = \frac{q^2 t^3}{(1-qt)^2} \cdot (B_1-A_1).
\end{equation}

\begin{figure}
\begin{center}
\includegraphics[width=114mm]{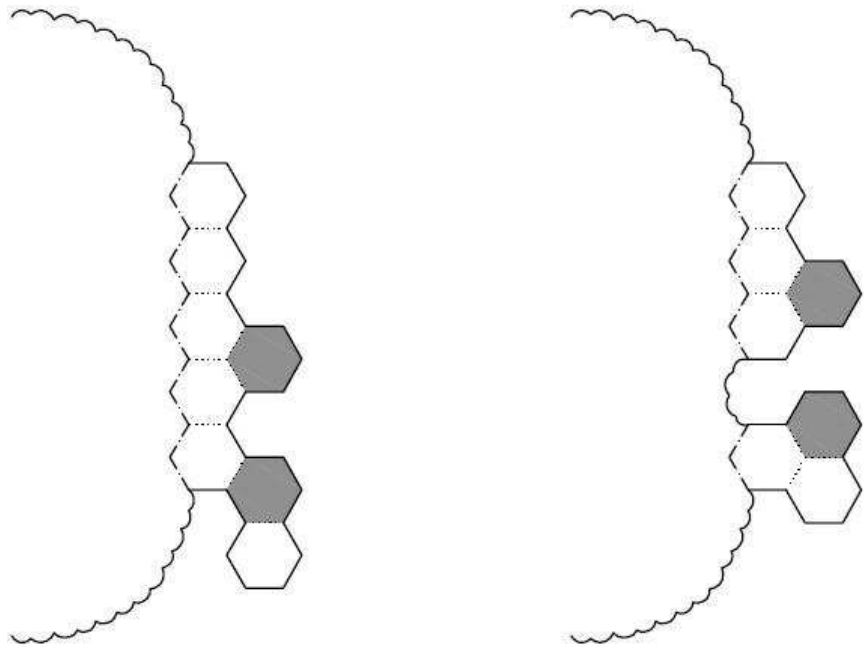}
\caption{The last two columns of two elements of $S_{\delta}$.}
\end{center}
\end{figure}

Now, let $P$ be an element of $S_{\epsilon}$. By the definition of $S_{\epsilon}$, $P$ is a polyomino with a one-part last column, but the body of $P$ is not a polyomino. Hence, in the second last column of $P$ there is a hole, and in the last column of $P$ there are two cells with which the hole is filled. In addition to this two-celled ``cork'', the last column contains $i \in \{0,\: 1,\: 2,\ldots \: \}$ cells lying below the ``cork'' and 
$j \in \{0,\: 1,\: 2,\ldots \: \}$ cells lying above the ``cork''. See Figure 8. Hence

\begin{equation}
A_{\epsilon}=\frac{1}{1-qt} \cdot q^2 t^2 \cdot \frac{1}{1-qt} \cdot C_1 = \frac{q^2 t^2}{(1-qt)^2} \cdot C_1.
\end{equation}

\begin{figure}
\begin{center}
\includegraphics[width=48.5mm]{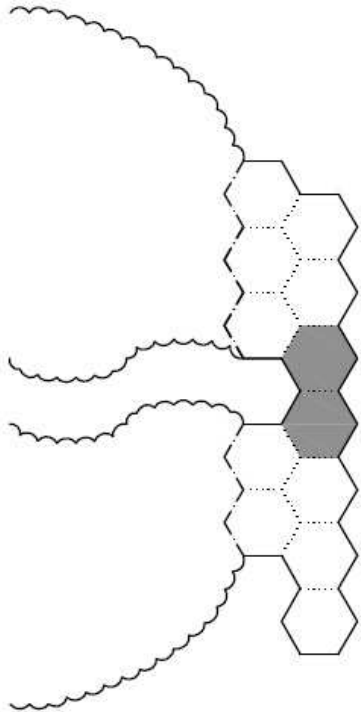}
\caption{The last two columns of an element of $S_{\epsilon}$.}
\end{center}
\end{figure}

If $P$ is an element of $S_{\zeta}$, then $P$ is a polyomino with a two-part last column, while the body of $P$ is not a polyomino. Once again, in the second last column of $P$ there is a hole, and in the last column of $P$ there are two cells with which the hole is filled. Let the lower component of the second last column consist of $i$ cells, and let the upper component of the second last column consist of $j$ cells. Now, if the two-celled cork belongs to the upper component of the last column, then it is impossible that $i-1$ or more cells of the last column lie between the cork and the hole of the last column. Namely, if $i-1$ or more cells were so situated, then the lower component of the last column would not be connected with the rest of $P$, and $P$ would not be a polyomino. See Figure 9.

\begin{figure}
\begin{center}
\includegraphics[width=128mm]{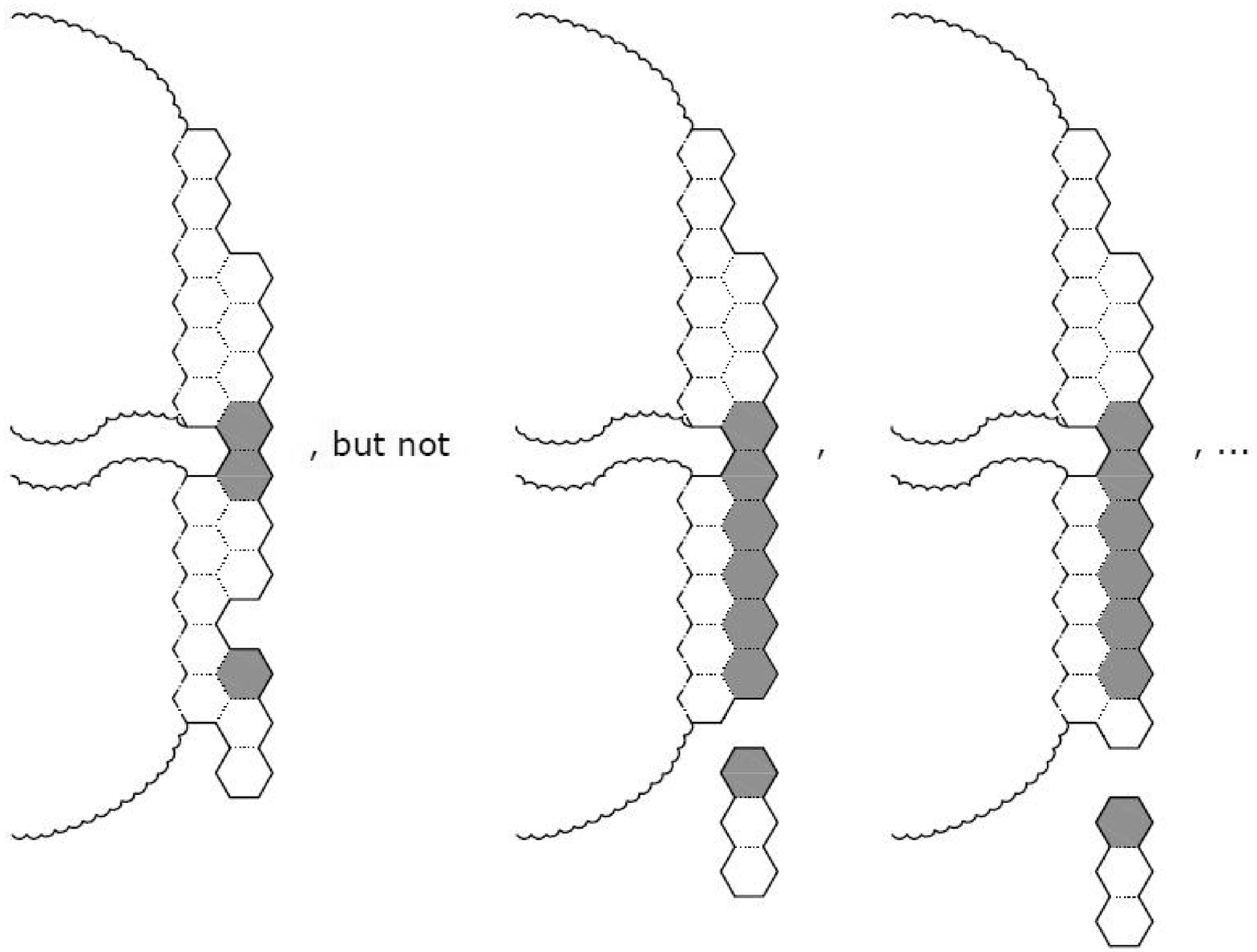}
\caption{The last two columns of an element of $S_{\zeta}$.}
\end{center}
\end{figure}

For a similar reason, if the two-celled cork belongs to the lower component of the last column, then it is impossible that $j-1$ or more cells of the last column lie between the cork and the hole of the last column. 

These remarks lead us to the following expression:

\begin{eqnarray*}
A_{\zeta} & = & \frac{q^3t^4}{(1-qt)^3} \cdot C_1 - \frac{q^3t^4}{(1-qt)^3} \cdot q^{-1}t^{-1}E(qt) \\
& & \mbox{} + \frac{q^3t^4}{(1-qt)^3} \cdot C_1 - \frac{q^3t^4}{(1-qt)^3} \cdot q^{-1}t^{-1}D(qt).
\end{eqnarray*}

When an element of $T$ is reflected about a horizontal axis, the area is preserved, whereas the height of the upper component of the last column becomes the height of the lower component of the last column. This means that $D(u)=E(u)$. We now have

\begin{equation}
A_{\zeta}=\frac{2q^3t^4}{(1-qt)^3} \cdot C_1 - \frac{2q^2t^3}{(1-qt)^3} \cdot D(qt).
\end{equation}

Since $A=A_{\alpha}+A_{\beta}+A_{\gamma}+A_{\delta}+A_{\epsilon}+A_{\zeta}$, equations (4.1)--(4.6) imply that

\begin{eqnarray}
A & = & \frac{qt}{1-qt} + \frac{qt}{(1-qt)^2} \cdot A_1 + \frac{qt}{1-qt} \cdot B_1 + \frac{q^2 t^3}{(1-qt)^2} \cdot (B_1-A_1) \nonumber \\
& & \mbox{} + \frac{q^2 t^2}{(1-qt)^2} \cdot C_1 + \frac{2q^3t^4}{(1-qt)^3} \cdot C_1 - \frac{2q^2t^3}{(1-qt)^3} \cdot D(qt).
\end{eqnarray}

Setting $t=1$, from equation (4.7) we get

\begin{eqnarray}
A_1 & = & \frac{q}{1-q} + \frac{q}{(1-q)^2} \cdot A_1 + \frac{q}{1-q} \cdot B_1 + \frac{q^2}{(1-q)^2} \cdot (B_1-A_1) \nonumber \\ 
& & \mbox{} + \frac{q^2}{(1-q)^2} \cdot C_1 + \frac{2q^3}{(1-q)^3} \cdot C_1 - \frac{2q^2}{(1-q)^3} \cdot D(q).
\end{eqnarray}

Differentiating equation (4.7) with respect to $t$ and then setting $t=1$, we get

\begin{eqnarray}
B_1 & = & \frac{q}{(1-q)^2} + \frac{q+q^2}{(1-q)^3} \cdot A_1 + \frac{q}{(1-q)^2} \cdot B_1 \nonumber \\
& & \mbox{} + \frac{3q^2-q^3}{(1-q)^3} \cdot (B_1-A_1) + \frac{2q^2}{(1-q)^3} \cdot C_1 + \frac{8q^3-2q^4}{(1-q)^4} \cdot C_1 \nonumber \\ 
& & \mbox{} - \frac{6q^2}{(1-q)^4} \cdot D(q) - \frac{2q^3}{(1-q)^3} \cdot D'(q).
\end{eqnarray}

\section{Setting up the functional equations for $C$, $D$ and $C_1$}

Now we turn to incomplete level one column-subconvex polyominoes. We have already observed that an incomplete level one column-subconvex polyomino always ends with a holed column.

The set $T_{\alpha}$ contains every two-part column (with one-celled hole) having $i \in \{1,\: 2,\: 3,\ldots \: \}$ cells below the hole and $j \in \{1,\: 2,\: 3,\ldots \: \}$ cells above the hole. Thus,

\begin{equation}
C_{\alpha} = \frac{qv}{1-qv} \cdot \frac{qu}{1-qu} = \frac{q^2uv}{(1-qu)(1-qv)}.
\end{equation}

If $P \in T_{\beta}$, then the body of $P$ lies in $S$. The hole of the last column has two possibilities: to coincide with the lower pivot cell of $P$ or to coincide with the upper pivot cell of $P$. Anyhow, the last column is made up of $i \in \{1,\: 2,\: 3,\ldots \: \}$ cells lying below the hole and $j \in \{1,\: 2,\: 3,\ldots \: \}$ cells lying above the hole. See Figure 10. Therefore,

\begin{equation}
C_{\beta} = 2 \cdot \frac{qv}{1-qv} \cdot \frac{qu}{1-qu} \cdot A_1 = \frac{2q^2uv}{(1-qu)(1-qv)} \cdot A_1.
\end{equation}

\begin{figure}
\begin{center}
\includegraphics[width=114mm]{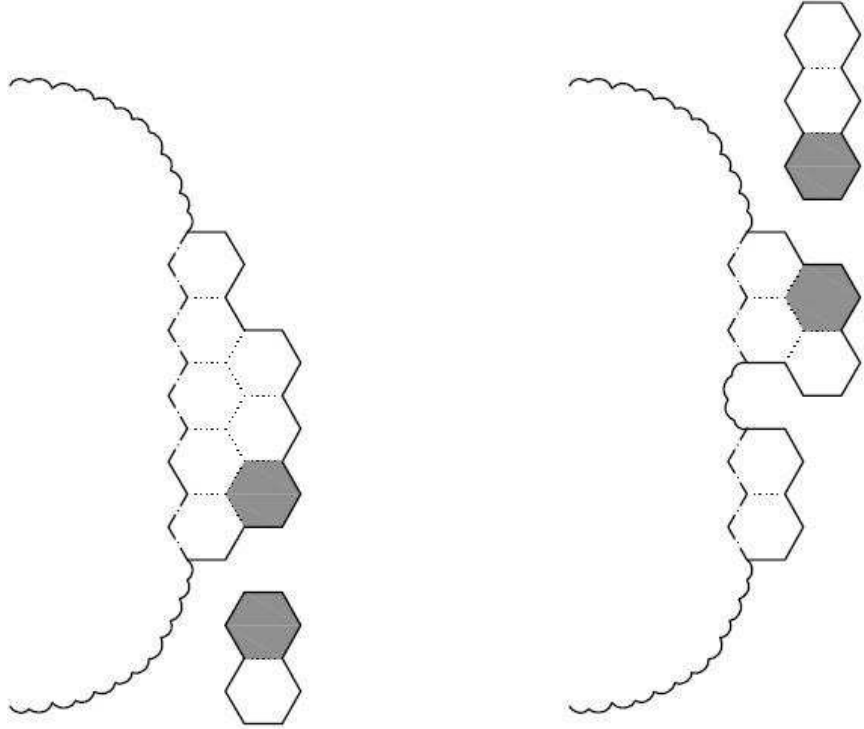}
\caption{The last two columns of two elements of $T_{\beta}$.}
\end{center}
\end{figure}

Now let $P \in T_{\gamma}$. The body of $P$ again lies in $S$. If the hole of the last column lies below the lower pivot cell of $P$, then the last column of $P$ is made up of:

\begin{itemize}
\item $i \in \{1,\: 2,\: 3,\ldots \: \}$ cells lying below the hole,
\item $j \in \{0,\: 1,\: 2,\ldots \: \}$ cells lying above the hole and below the lower pivot cell,
\item the lower pivot cell, and
\item $k \in \{0,\: 1,\: 2,\ldots \: \}$ cells lying above the lower pivot cell.
\end{itemize}

If the hole of the last column lies above the upper pivot cell of $P$, then the last column of $P$ is made up of:

\begin{itemize}
\item $i \in \{1,\: 2,\: 3,\ldots \: \}$ cells lying above the hole,
\item $j \in \{0,\: 1,\: 2,\ldots \: \}$ cells lying below the hole and above the upper pivot cell,
\item the upper pivot cell, and
\item $k \in \{0,\: 1,\: 2,\ldots \: \}$ cells lying below the upper pivot cell.
\end{itemize}

See Figure 11. Altogether,

\begin{equation}
C_{\gamma} = \frac{q^2uv}{(1-qu)^2(1-qv)} \cdot A_1 + \frac{q^2uv}{(1-qu)(1-qv)^2} \cdot A_1.
\end{equation}

\begin{figure}
\begin{center}
\includegraphics[width=114mm]{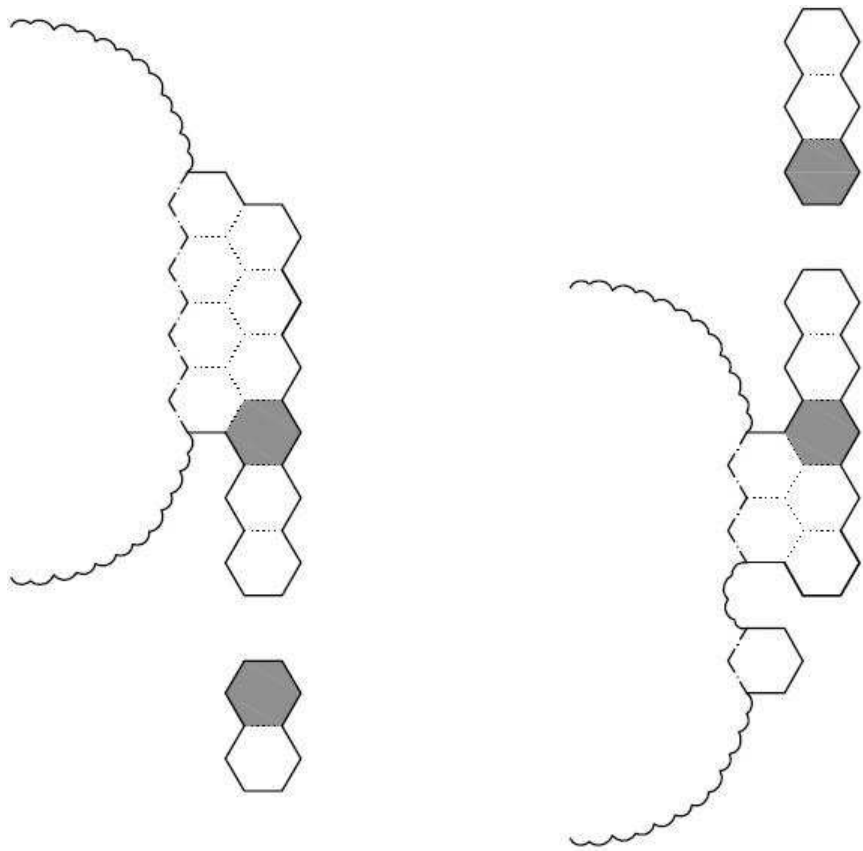}
\caption{The last two columns of two elements of $T_{\gamma}$.}
\end{center}
\end{figure}

If $P \in T_{\delta}$, then the body of $P$ lies in $T$. The second last and last columns of $P$ both have a hole. The hole of the last column is either the lower right neighbour or the upper right neighbour of the hole of the second last column. In the last column, there are $i \in \{1,\: 2,\: 3,\ldots \: \}$ cells below the hole and $j \in \{1,\: 2,\: 3,\ldots \: \}$ cells above the hole. See Figure 12. Hence,

\begin{equation}
C_{\delta} = 2 \cdot \frac{qv}{1-qv} \cdot \frac{qu}{1-qu} \cdot C_1 = \frac{2q^2uv}{(1-qu)(1-qv)} \cdot C_1.
\end{equation}

\begin{figure}
\begin{center}
\includegraphics[width=\textwidth]{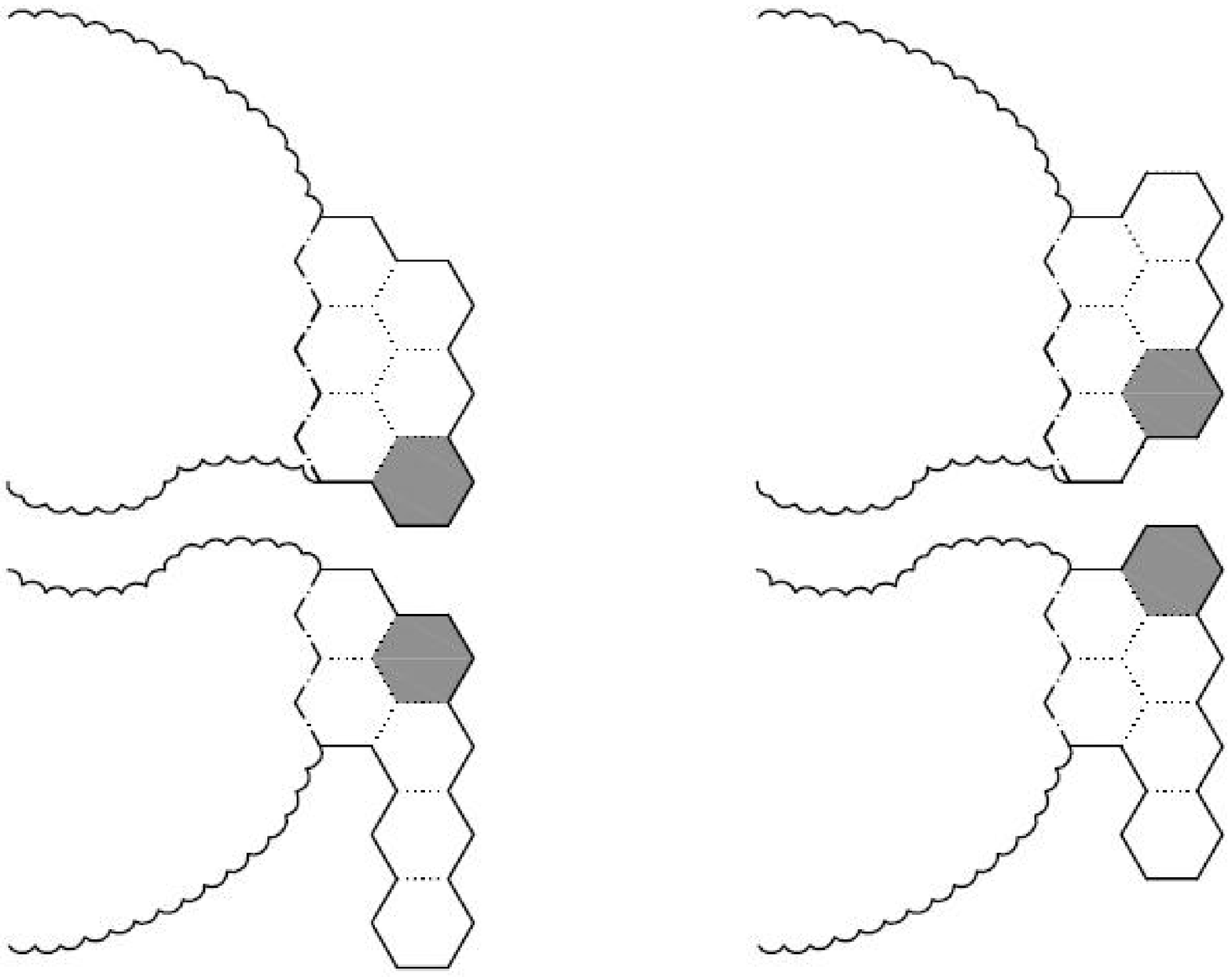}
\caption{The last two columns of two elements of $T_{\delta}$.}
\end{center}
\end{figure}

Let $P \in T_{\epsilon}$. Once again, the second last and last columns of $P$ both have a hole. However, to the right of the hole of the second last column, there are two cells which both belong to $P$. If this two-celled cork is contained in the lower component of the last column, and if the upper component of the second last column consists of $j$ cells, then it is necessary that at least $j-1$ cells of the last column lie above the cork and below the hole of the last column. Otherwise the upper component of the last column would be connected with the rest of $P$, and $P$ would be a polyomino. (That cannot happen because $P$ is an element of the set $T$, and the elements of $T$ are not polyominoes.) See Figure 13.

\begin{figure}
\begin{center}
\includegraphics[width=\textwidth]{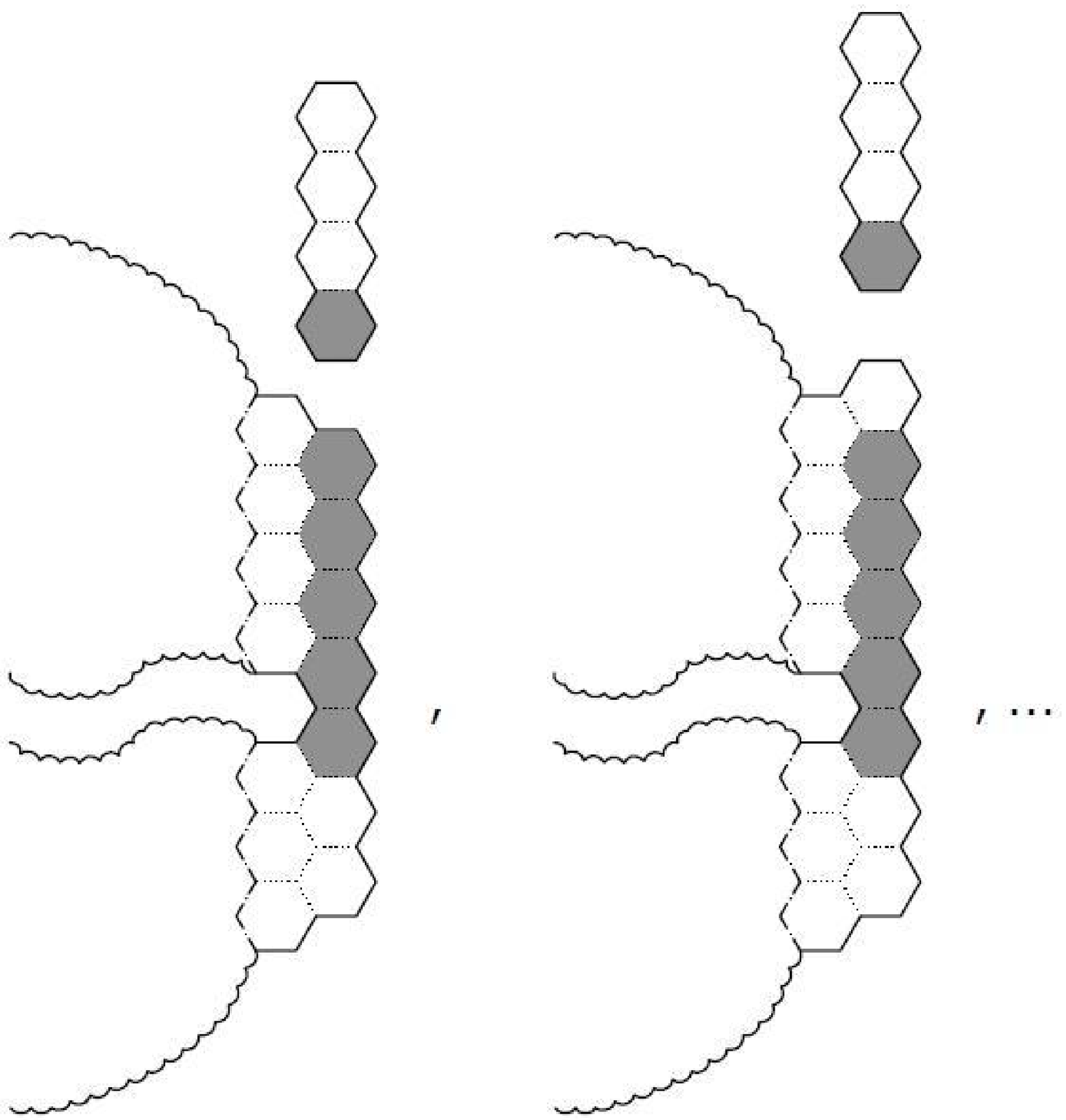}
\caption{The last two columns of the elements of $T_{\epsilon}$.}
\end{center}
\end{figure}

Similarly, if the two-celled cork is contained in the upper component of the last column, and if the lower component of the second last column consists of $i$ cells, then it is necessary that at least $i-1$ cells of the last column lie below the cork and above the hole of the last column. If it were not so, then $P$ would be a polyomino.

These remarks made, we conclude that

\begin{eqnarray}
C_{\epsilon} & = & \frac{q^3uv^2}{(1-qu)(1-qv)^2} \cdot q^{-1}v^{-1}D(qv) + \frac{q^3u^2v}{(1-qu)^2(1-qv)} \cdot q^{-1}u^{-1}E(qu) \nonumber \\
& = & \frac{q^2uv}{(1-qu)(1-qv)^2} \cdot D(qv) + \frac{q^2uv}{(1-qu)^2(1-qv)} \cdot D(qu).
\end{eqnarray}

Since $C=C_{\alpha}+C_{\beta}+C_{\gamma}+C_{\delta}+C_{\epsilon}$, equations (5.1)--(5.5) imply that

\begin{eqnarray}
C & = & \frac{q^2uv}{(1-qu)(1-qv)} + \frac{2q^2uv}{(1-qu)(1-qv)} \cdot A_1 + \frac{q^2uv}{(1-qu)^2(1-qv)} \cdot A_1 \nonumber \\
& & \mbox{} + \frac{q^2uv}{(1-qu)(1-qv)^2} \cdot A_1 + \frac{2q^2uv}{(1-qu)(1-qv)} \cdot C_1 \nonumber \\
& & \mbox{} + \frac{q^2uv}{(1-qu)(1-qv)^2} \cdot D(qv) + \frac{q^2uv}{(1-qu)^2(1-qv)} \cdot D(qu).
\end{eqnarray}

Setting $v=1$, from equation (5.6) we get

\begin{eqnarray}
D(u) & = & \frac{q^2u}{(1-q)(1-qu)} + \frac{q^2u}{(1-q)(1-qu)^2} \cdot A_1 \nonumber \\
& & \mbox{} + \frac{(3-2q)q^2u}{(1-q)^2(1-qu)} \cdot A_1 + \frac{2q^2u}{(1-q)(1-qu)} \cdot C_1 \nonumber \\ 
& & \mbox{} + \frac{q^2u}{(1-q)^2(1-qu)} \cdot D(q) + \frac{q^2u}{(1-q)(1-qu)^2} \cdot D(qu).
\end{eqnarray}

We have $D(1)=C(q,1,1)=C_1$. So, when we set $u=1$, equation (5.7) turns into

\begin{equation}
C_1=\frac{q^2}{(1-q)^2} + \frac{4q^2-2q^3}{(1-q)^3} \cdot A_1 + \frac{2q^2}{(1-q)^2} \cdot C_1 + \frac{2q^2}{(1-q)^3} \cdot D(q).
\end{equation}

\section{Solving the functional equations}

For convenience, we firstly define an extra series $F$. The definition is

\begin{equation}
F=1 + \frac{3-2q}{1-q} \cdot A_1 + 2C_1 + \frac{1}{1-q} \cdot D(q).
\end{equation}

Now equation (5.7) can be written as

\begin{equation}
D(u)=\frac{q^2u}{(1-q)(1-qu)^2} \cdot A_1 + \frac{q^2u}{(1-q)(1-qu)} \cdot F + \frac{q^2u}{(1-q)(1-qu)^2} \cdot D(qu).
\end{equation}

The next step of the upgraded Temperley method is to iteratively remove the $D(qu)$ term from the right-hand side of equation (6.2). Namely, substituting $qu$ for $u$ in equation (6.2) gives

\begin{displaymath}
D(qu)=\frac{q^3u}{(1-q)(1-q^2u)^2} \cdot A_1 + \frac{q^3u}{(1-q)(1-q^2u)} \cdot F + \frac{q^3u}{(1-q)(1-q^2u)^2} \cdot D(q^2u).
\end{displaymath}

Replacing $D(qu)$ of equation (6.2) by this latter expression, followed by a bit of rearranging, results in

\begin{eqnarray*}
D(u) & = & \left[\frac{q^2u}{(1-q)(1-qu)^2} + \frac{q^{2+3}u^2}{(1-q)^2(1-qu)^2(1-q^2u)^2} \right] \cdot A_1 \\
& & \mbox{} + \left[\frac{q^2u}{(1-q)(1-qu)} + \frac{q^{2+3}u^2}{(1-q)^2(1-qu)^2(1-q^2u)} \right] \cdot F \\
& & \mbox{} + \frac{q^{2+3}u^2}{(1-q)^2(1-qu)^2(1-q^2u)^2} \cdot D(q^2u).
\end{eqnarray*}

After the next iteration, in each of the square brackets there is a sum of three terms, and the argument of the final $D$ is $q^3u$ instead of $q^2u$. After infinitely many iterations, we have

\begin{eqnarray}
D(u) & = & \left\{ \sum_{i=1}^{\infty} \frac{q^{\frac{i(i+3)}{2}}u^i}{(1-q)^i \cdot \left[ \prod_{k=1}^{i}(1-q^ku) \right]^2} \right\} \cdot A_1 \nonumber \\
& & \mbox{} + \left\{ \sum_{i=1}^{\infty} \frac{q^{\frac{i(i+3)}{2}}u^i}{(1-q)^i \cdot \left[ \prod_{k=1}^{i-1}(1-q^ku) \right]^2 \cdot (1-q^iu)} \right\} \cdot F.
\end{eqnarray}

The right-hand side of equation (6.3) involves no $D$ because $lim_{n \rightarrow \infty} D(q^nu)=0$. The reason why this limit is zero is that the lowest power of $q$ occurring in $D(q^nu)$ is $n+2$.

Setting $u=q$, from equation (6.3) we get

\begin{eqnarray}
D(q) & = & \left\{ \sum_{i=1}^{\infty} \frac{q^{\frac{i(i+5)}{2}}}{(1-q)^i \cdot \left[ \prod_{k=1}^{i}(1-q^{k+1}) \right]^2} \right\} \cdot A_1 \nonumber \\ 
& & \mbox{} + \left\{ \sum_{i=1}^{\infty} \frac{q^{\frac{i(i+5)}{2}}}{(1-q)^i \cdot \left[ \prod_{k=1}^{i-1}(1-q^{k+1}) \right]^2 \cdot (1-q^{i+1})} \right\} \cdot F.
\end{eqnarray}

``Logarithmically'' differentiating\footnote{By logarithmic differentiation we mean the use of the formula $\varphi '=\varphi \cdot [ln(\varphi)]'$.} equation (6.3) with respect to $u$ and then setting $u=q$, we obtain

\begin{eqnarray}
D'(q) & = & \left\{ \sum_{i=1}^{\infty} \frac{q^{\frac{i(i+5)}{2}}}{(1-q)^i \cdot \left[ \prod_{k=1}^{i}(1-q^{k+1}) \right]^2} \cdot \left(\frac{i}{q} +2 \cdot \sum_{j=1}^i \frac{q^j}{1-q^{j+1}} \right) \right\} \cdot A_1 \nonumber \\ 
& & \mbox{} + \left\{ \sum_{i=1}^{\infty} \frac{q^{\frac{i(i+5)}{2}}}{(1-q)^i \cdot \left[ \prod_{k=1}^{i-1}(1-q^{k+1}) \right]^2 \cdot (1-q^{i+1})} \right. \nonumber \\
& & \left. \cdot \left( \frac{i}{q} +2 \cdot \sum_{j=1}^{i-1} \frac{q^j}{1-q^{j+1}} + \frac{q^i}{1-q^{i+1}} \right) \right\} \cdot F.
\end{eqnarray}

Equations (4.8), (4.9), (5.8), (6.1), (6.4) and (6.5) make up a system of six \textit{linear equations} in six unknowns: $A_1$, $B_1$, $C_1$, $D(q)$, $D'(q)$ and $F$. That linear system was readily solved by the computer algebra package \textit{Maple}. Of course, the most interesting component of the solution is $A_1$, the area generating function for level one column-subconvex polyominoes. We state the formula for $A_1$ as a theorem.

\begin{theo} The area generating function for level one column-subconvex polyominoes is given by

\begin{displaymath}
A_1=\frac{\sum_{n=1}^{3}num_n}{\sum_{n=1}^{6}den_n},
\end{displaymath}

\noindent where

\begin{eqnarray*}
num_1 & = & q-8q^{2}+24q^{3}-32q^{4}+17q^{5}+4q^{6}-8q^{7}+2q^{8}, \\
num_2 & = & (-q+5q^{2}-13q^{3}+23q^{4}-22q^{5}+12q^{6}-2q^{7}) \cdot \beta , \\
num_3 & = & (-2q^{4}+8q^{5}-12q^{6}+8q^{7}-2q^{8}) \cdot \delta , \\
& & \\
den_1 & = & 1-11q+46q^{2}-93q^{3}+88q^{4}-27q^{5}-24q^{6}+19q^{7}-3q^{8}, \\
den_2 & = & (2q^{2}-8q^{3}+8q^{4}-4q^{5}-6q^{6}+4q^{7}) \cdot \alpha , \\
den_3 & = & (-1+10q-34q^{2}+67q^{3}-81q^{4}+54q^{5}-16q^{6}+q^{7}) \cdot \beta , \\
den_4 & = & (2q^{4}-8q^{5}+8q^{6}-2q^{8}) \cdot \gamma , \\
den_5 & = & (6q^{4}-22q^{5}+34q^{6}-22q^{7}+4q^{8}) \cdot \delta , \\
den_6 & = & (2q^{4}-6q^{5}+10q^{6}-6q^{7}) \cdot (\alpha\delta-\beta\gamma), 
\end{eqnarray*}

\begin{eqnarray*}
\alpha & = & \sum_{i=1}^{\infty} \frac{q^{\frac{i(i+5)}{2}}}{(1-q)^{i}\left[ \prod_{k=1}^{i}(1-q^{k+1})\right]^{2}} \ , \\
& & \\
\beta & = & \sum_{i=1}^{\infty} \frac{q^{\frac{i(i+5)}{2}}}{(1-q)^{i}\left[ \prod_{k=1}^{i-1}(1-q^{k+1})\right]^{2}(1-q^{i+1})} \ , \\
& & \\
\gamma & = & \sum_{i=1}^{\infty} \frac{q^{\frac{i(i+5)}{2}} \left(\frac{i}{q}+2 \sum_{j=1}^{i}\frac{q^j}{1-q^{j+1}} \right) }{(1-q)^{i}\left[ \prod_{k=1}^{i}(1-q^{k+1})\right]^{2}} \ , \\
& & \\
\delta & = & \sum_{i=1}^{\infty} \frac{q^{\frac{i(i+5)}{2}} \left(\frac{i}{q}+2 \sum_{j=1}^{i-1}\frac{q^j}{1-q^{j+1}}+\frac{q^{i}}{1-q^{i+1}}\right) }{(1-q)^{i}\left[ \prod_{k=1}^{i-1}(1-q^{k+1})\right]^{2}(1-q^{i+1})} \ .
\end{eqnarray*}
\end{theo}

From the formula just stated, one easily finds that

\begin{eqnarray*}
A_1 & = & q+3q^2+11q^3+44q^4+184q^5+786q^6+3391q^7+14683q^8 \\
& & +63619q^9+275506q^{10}+1192134q^{11}+5154794q^{12}+\ldots \ .
\end{eqnarray*}

We expanded $A_1$ in a Taylor series to 250 terms, and analysed the series by the method of differential approximants \cite{critical} using second-order approximants, that is to say, approximants given by solutions of inhomogeneous second degree ordinary differential equations. From this analysis, we found that the dominant singularity of $A_1$ is a simple pole, located at $q=q_c=0.2315276132 \: $. (Note that we only needed some 20 series terms to establish this---the additional terms merely provided higher accuracy and confirmation of our initial analysis). We refined this estimate by using \textit{Maple} to locate the position of the denominator zero. That is, by expanding the series $\sum_{n=1}^6 den_n$ to more and more terms, more and more accurate numerical solutions of $\sum_{n=1}^6 den_n = 0$ were obtained. In this way, we found $q_c$ to be $0.231527613159$. We could obtain much higher accuracy if necessary. It is likely that this number is algebraic (such is usually the case with exact solutions), but we have been unable to conjecture its exact value. Note also that the numerator is positive for $0 < q < 0.8,$ so there is no possibility that this denominator zero cancels with the numerator.

The growth constant is the reciprocal of the dominant singularity. Thus, the growth constant of level one column-subconvex polyominoes is about $4.319139$. For comparison, the growth constant of column-convex polyominoes is $3.863131$, the growth constant of level one cheesy polyominoes is $4.114908$, and the growth constant of all polyominoes is $5.183148$. (The latter two growth constants were found in \cite{semi} and \cite{Voege}, respectively.) From the result that $\tau = \sup_{n \ge 1} a_n^{1/n},$ we also have the quite good lower bound (based on 250 terms) $\tau >  4.283006.$ We can also calculate the amplitude, so writing the generating function as $A_1 = \sum_n a_nq^n,$ then $a_n \sim c_1 \cdot \tau^n,$ we can estimate the amplitude $c_1$ from the sequence of quotients $a_n/\tau^n.$ In this way we estimate $c_1 = 0.1224281005.$

\section{Level two column-subconvex polyominoes}

In the just-finished enumeration of level one column-subconvex polyominoes, we considered altogether $11$ cases. Namely, we partitioned the set $S$ into $6$ subsets and the set $T$ into $5$ subsets; $6+5$ equals $11$. 

Let us write $A_2$ to denote the area generating function for level two column-subconvex polyominoes. We found a formula for $A_2$ as well, but that goal was achieved through considering as much as $37$ cases. We had to struggle against a number of complicated expressions, and in the end we had to solve a system of $16$ linear equations in $16$ unknowns. (For comparison, the computation of $A_1$ was completed by solving a system of $6$ linear equations in $6$ unknowns.) Consequently, the formula for $A_2$ is much bulkier than the formula for $A_1$. To save this journal's space, we have chosen to state the formula for $A_2$ in the electronic form only \cite{worksheet}.

However, the Taylor series expansion of $A_2$ is

\begin{eqnarray*}
A_2 & = & q+3q^2+11q^3+44q^4+186q^5+812q^6+3614q^7+16254q^8 \\
& & \mbox{}+73464q^9+332603q^{10}+1505877q^{11}+6813301q^{12}+\ldots
\end{eqnarray*}

\noindent and the critical point of $A_2$ is at $q=q_c=0.221755050048$. This was obtained in the same way as described above for level one column-subconvex polyominoes, but based on a series of 153 terms. Thus the growth constant of level two column-subconvex polyominoes is about $4.509480\: $. For comparison, the growth constant of level two cheesy polyominoes is $4.231836$ \cite{semi}. As above, we can also give the rigorous bound $ \tau > 4.441222.$ We can also estimate the amplitude $c_2 = 0.0969488405,$ so that the coefficient of the $n^{th}$ term of the generating function $A_2 = \sum_n a_nq^n$ is $a_n \sim c_2 \cdot \tau^n.$ 

As stated in \cite{Voege}, the area generating function for all polyominoes is

\begin{eqnarray*}
& & q+3q^2+11q^3+44q^4+186q^5+814q^6+3652q^7+16689q^8 \\
& & \mbox{}+77359q^9+362671q^{10}+1716033q^{11}+8182213q^{12}+\ldots \ .
\end{eqnarray*}

Indeed, a quick drawing confirms that a polyomino must have at least $5$ (resp. $6$) cells in order not to be a level one (resp. two) column-subconvex polyomino. See Figure 14.

\begin{figure}
\begin{center}
\includegraphics[width=\textwidth]{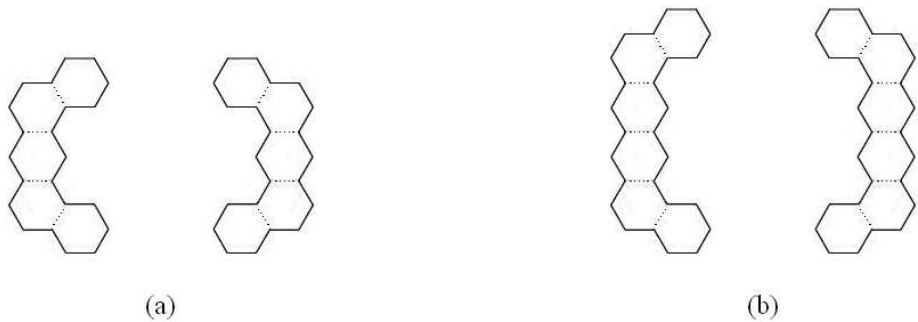}
\caption{(a) The two smallest instances of a polyomino which is not a level one column-subconvex polyomino. (b) The two smallest instances of a polyomino which is not a level two column-subconvex polyomino.}
\end{center}
\end{figure}

We have not tried to enumerate level three column-subconvex polyominoes. Our non-rigorous estimate is that, in order to enumerate this latter model by area, one would have to consider at least $80$ cases.

\section{Conclusion}

We have defined a class of polyominoes that interpolates between column-convex polyominoes and all polyominoes. The former have been solved, while the latter remain unsolved. For now, our interpolating class (we call it level $m$ column-subconvex polyominoes) remains solved up to a certain point. Namely, we have solved the cases $m=1$ and $m=2$. Column-convex polyominoes correspond to the case $m=0$. In both cases $m=1,\: 2$, the generating function has a simple pole singularity, located at $q=q_c=0.2315\ldots$ and $0.2217\ldots$ respectively. For all polyominoes, the corresponding singularity is at $q=q_c(\mathrm{polyomino})=0.192932\ldots \:$, and the singularity is of the form $const. \cdot \vert \mathrm{log}(q_c-q) \vert$, rather than a simple pole \cite{Voege}. For all finite values of $m$ we expect the generating function of level $m$ column-subconvex polyominoes to have a simple pole, while the singularity position is expected to be a monotone decreasing function of $m$, with a limiting value as $m$ tends to infinity of $q^{*}>q_c(\mathrm{polyomino})$. We have also given the rigorous lower bounds $\tau > 4.283006$ and $\tau > 4.441222$ for the growth constants of level~1 and level~2 column-subconvex polyominoes respectively.

\end{document}